\documentstyle{article}

\begin{document}

\begin{center}
\begin{flushleft}
{\Large \bf A New Aspect of Representations of $U_q(\hat{sl_2})$\\
---Generic Case}\\
\vspace{4mm}
Xufeng Liu\\
\vspace{3mm}
Department of Mathematics, The University of Melbourne\\ Parkville, Victoria 3052,
Australia\\
\vspace{3mm}
Department of Mathematics, Peking University\\ Beijing 100871, P.R.China
\footnote{Permanent mailing address}
\end{flushleft}
\end{center}

\vspace{20mm} 
\begin{flushleft}
Abstract.An identity is derived from the tensor product representation $V_{m}(x)\otimes
V_{n}(y)$ of $U_q(\hat{sl_2})$ and a new basis of $V_{m}(x)\otimes V_{n}(y)$
is established.
\end{flushleft}
\newpage

In this paper we take $q$ to be a complex number instead of a formal
variable, unless otherwise pointed out explicitly.For any complex number $q$
and any integers $r$ and $s$ we use the notations 
\begin{eqnarray*}
[r]_q&=&\frac{q^{r}-q^{-r}}{q-q^{-1}}, \\
\left [r\right ]_{q}!&=&\prod_{s=1}^{r} [s]_{q}, \\
\left [\begin{array}{c}r \\ s\end{array}
\right ]_{q}&=&\frac{[r]_{q}!}{[s]_{q}![r-s]_{q}!},r\geq s.
\end{eqnarray*}
and we set $[0]_{q}!=1$.We also assume that $q$ is not a root of unity.

\vspace{4mm}

{\bf 1.Some Basic Facts}\\
\vspace{1mm}

In this section we list some basic facts about $U_q(\hat{sl_2})$ and fix
notations.For details we refer the readers to Ref.[1].

 {\bf Definition 1.1.} The quantum affine algebra $U_q(\hat{sl_2})$  is the associative algebra over $C$
with generators $e_{i},f_{i},K_{i}$ and $K_{i}^{-1}$($i=0,1$) and the
following relations: 
\begin{eqnarray*}
&&K_{i}K_{i}^{-1} = K_{i}^{-1}K_{i}=1, \\
&&K_0K_1 = K_1K_0, \\
&&K_{i}e_{i}K_{i}^{-1}=q^2e_{i},K_{i}f_{i}K_{i}^{-1}=q^{-2}f_{i}, \\
&&K_{i}e_{j}K_{i}^{-1}=q^{-2}e_{j},K_{i}f_{j}K_{i}^{-1}=q^2f_{j},i\ne j, \\
&&\left [e_{i},f_{i}\right ]=\frac{K_{i}-K_{i}^{-1}}{q-q^{-1}}, \\
&&\left [e_0,f_1\right ]=\left [e_1,f_0\right ]=0,
\end{eqnarray*}
\begin{eqnarray*}
e_{i}^3e_{j}-[3]_qe_{i}^2e_{j}e_{i}+[3]_qe_{i}e_{j}e_{i}^2-e_{j}e_{i}^3&=&0,
\\
f_{i}^3f_{j}-[3]_qf_{i}^2f_{j}f_{i}+[3]_qf_{i}f_{j}f_{i}^2-f_{j}f_{i}^3&=&0.(i\ne j).
\end{eqnarray*}
Moreover, $U_q(\hat{sl_2})$ is a Hopf algebra over $C$ with the
comultiplication 
\begin{eqnarray*}
&&\triangle (e_{i})=e_{i}\otimes K_{i}+1\otimes e_{i}, \\
&&\triangle (f_{i})=f_{i}\otimes 1+K_{i}^{-1}\otimes f_{i}, \\
&&\triangle (K_{i})=K_{i}\otimes K_{i},\triangle
(K_{i}^{-1})=K_{i}^{-1}\otimes K_{i}^{-1}, \\
\end{eqnarray*}
and the antipode 
\begin{eqnarray*}
&&S(K_{i})=K_{i}^{-1},S(K_{i}^{-1})=K_{i}, \\
&&S(e_{i})=-e_{i}K_{i}^{-1},S(f_{i})=-K_{i}f_{i}.
\end{eqnarray*}

{\bf Definition 1.2.} The quantum algebra $U_q(sl_2)$ is the associative algebra
over $C$ with generators $e,f,K$ and $K^{-1}$ and the following relations: 
\begin{eqnarray*}
&&KK^{-1}=K^{-1}K=1, \\
&&KeK^{-1}=q^2e,KfK^{-1}=q^{-2}f, \\
&&\left [e,f\right ]=\frac{K-K^{-1}}{q-q^{-1}}.
\end{eqnarray*}
It is a Hopf algebra over $C$ with the following comultiplication $\triangle$
and antipode $S$: 
\begin{eqnarray*}
&&\triangle e=e\otimes K+1\otimes e, \\
&&\triangle f=f\otimes 1+K^{-1}\otimes f, \\
&&\triangle K=K\otimes K,\triangle K^{-1}=K^{-1}\otimes K^{-1}, \\
&&S(K)=K^{-1},S(K^{-1})=K,S(e)=-eK^{-1},S(f)=-Kf.
\end{eqnarray*}

There is an associative algebra homomorphism, known as evaluation map, from $%
U_q(\hat{sl_2})$ to $U_q(sl_2)$[2].

{\bf Definition 1.3.}For any $x\in
C\backslash\{0\}$, the evaluation map $ev_x$  from $U_q(\hat{sl_2})$ to $%
U_q(sl_2)$ is the associative algebra homomorphism such that 
\begin{eqnarray*}
&&ev_x(e_0)=q^{-1}xf, ev_x(e_1)=e, \\
&&ev_x(f_0)=qx^{-1}e,ev_x(f_1)=f, \\
&&ev_x(K_0)=K^{-1},ev_x(K_1)=K.
\end{eqnarray*}

It is clear that modules of $U_q(\hat{sl_2})$ can be obtained by pulling
back modules of $U_q(sl_2)$ by the homomorphism $ev_x$.We denote by $V(x)$
the pull back module of $U_q(\hat{sl_2})$ of a module $V$ of $U_q(sl_2)$ by
the evaluation map $ev_x$.

{\bf Definition 1.4.} Let $V$ and $W$ be two modules of $U_q(sl_2)$.The $U_q(\hat{%
sl_2})$ module $W(y)\otimes V(x)$ is the vector space $W\otimes V$ with the
module structure defined through the following action 
\[
g(w\otimes v)\stackrel{\triangle}{=}(ev_y\otimes ev_x)\triangle g(w\otimes
v),\forall g\in U_q(\hat{sl_2}).
\]

We consider the module $V_m(x)\otimes V_n(y)$ of $U_q(\hat{sl_2})$, where $%
V_n$ is the standard $n+1$ dimensional module of $U_q(sl_2)$. There 
is a basis $\{v_i|i=0,1,\cdots ,m\}$ of $V_m$ such that its module structure
is  defined through the following actions: 
\begin{eqnarray*}
Kv_i &=&q^{m-2i}v_i, \\
fv_i &=&[i+1]_qv_{i+1}, \\
ev_i &=&[m+1-i]_qv_{i-1},i=0,1,\cdots ,m.
\end{eqnarray*}
Here we have used the convention $v_{-1}=v_{m+1}=0$. From now on we denote
by $\{w_i|i=0,1,\cdots ,n\}$ the basis of $V_n$ satisfying these equations
to distinguish it from that of $V_m$.Finally, we recall from [2] that, as a
representation of $U_q(sl_2)$, 
\[
V_m\otimes V_n\cong V_{m+n}\oplus V_{m+n-2}\oplus \cdots \oplus V_{|m-n|}.
\]
\vspace{2mm}
{\bf 2.An Identity}
\vspace{1mm}

In this section we will derive an identity from the representation \\$%
V_{m}(x)\otimes V_{n}(y)$.

For any integer $l\leq min\{m,n\}$ ,let 
\[
\Omega _l=\sum_{i=0}^lc_{i,l-i}v_i\otimes w_{l-i},\\ \Phi
_l=\sum_{i=0}^ld_{i,l}v_{m-l+i}\otimes w_{n-i},
\]
where 
\begin{eqnarray*}
c_{i,l-i} &=&(-1)^iq^{i(2l-n-i-1)}\prod_{j=0}^i\frac{[n-l+j]_q}{[m-j+1]_q},
\\
d_{i,l} &=&(-1)^iq^{i(-m+2l-i-1)}\prod_{j=0}^i\frac{[m-l+j]_q}{[n-j+1]_q}.
\end{eqnarray*}
It is readily verified that 
\[
e\Omega _l=0,f\Phi _l=0.
\]
Besides, $\Omega _l$ are the only highest vectors (up to a scalar multiple)
in $V_m(x)\otimes V_n(y)$ with respect to $U_q({sl_2})$ and $\Phi _l$ the
only lowest vectors. As a matter of fact, $\Omega _l$ generates the $%
U_q(sl_2)$ submodule $V_{m+n-2l}$. It then follows that $f^{m+n-2l}\Omega _l$
must be a scalar multiple of $\Phi _l$: 
\[
f^{m+n-2l}\Omega _l=\alpha _l\Phi _l,
\]
where $\alpha _l$ is a constant.

We find that there are two ways, direct and indirect ones, of determining $%
\alpha _{l}$. Thus by equating the two results from different ways we are
able to establish an identity.

First, let us take the direct way. For any positive integer $k$ we have 
\[
\triangle f^{k}=\sum_{j=0}^{k} q^{-j(k-j)}\left [%
\begin{array}{c}
k \\ 
j
\end{array}
\right ]_{q}K^{-j}f^{k-j}\otimes f^{j}. 
\]
Using this formula, after some elementary calculation, we get the following
coefficient of the term $v_{m-l}\otimes w_{n}$ in $f^{m+n-2l}\Omega _{l}$: 
\begin{eqnarray*}
&&q^{-l(n-l)}\sum_{i=0}^{min\{l,m-l\}} (-1)^{i}q^{-i}\frac
{[m+n-2l]_{q}![m-l]_{q}![n]_{q}!}
{[m-l-i]_{q}![n-l+i]_{q}![i]_{q}![l-i]_{q}!}\times \\
&&\prod_{j=0}^{i} \frac{[n-l+j]_{q}}{[m-j+1]_{q}}.
\end{eqnarray*}
On the other hand, the coefficient of the same term in $\Phi _{l} $ is 
\[
d_{0,l}=\frac{[m-l]_{q}}{[n+1]_{q}}. 
\]
It follows that 
\begin{eqnarray*}
\alpha _{l}&=& q^{-l(n-l)}\frac{[n+1]_{q}}{[m-l]_{q}} \sum_{i=0}^{min\{l,m-l%
\}} (-1)^{i}q^{-i}\frac {[m+n-2l]_{q}![m-l]_{q}![n]_{q}!}
{[m-l-i]_{q}![n-l+i]_{q}![i]_{q}![l-i]_{q}!}\times \\
&&\prod_{j=0}^{i} \frac{[n-l+j]_{q}}{[m-j+1]_{q}}.
\end{eqnarray*}

Now let us turn to probe the detour. We need the following results which can
be proved by direct calculation.

{\bf Lemma 2.1.}For $l\geq 1$ 
\[
e^2e_0\Omega _l=[2]_q[n-l]_qq^{-1}(xq^m-yq^{-n+2l-2})\Omega _{l-1}.
\]

{\bf Lemma 2.2.}For $l\geq 1$ 
\[
e_0\Phi _l=[m-l]_qq^{-1}(xq^n-yq^{-m+2l-2})\Phi _{l-1}.
\]

It follows from Lemma 2.1 that 
\[
e_0\Omega _l=c_{l-1}f^2\Omega _{l-1}+c_lf\Omega _l+c_{l+1}\Omega _{l+1},
\]
where 
\[
c_{l-1}=\frac{[n-l]_qq^{-1}(xq^m-yq^{-n+2l-2})}{[m+n-2l+1]_q[m+n-2l+2]_q},
\]
and $c_l$ and $c_{l+1}$ are two other constants.We then have 
\[
e_0f^{m+n-2l}\Omega _l=c_{l-1}f^{m+n-2l+2}\Omega _{l-1},
\]
namely, 
\[
\alpha _le_0\Phi _l=c_{l-1}\alpha _{l-1}\Phi _{l-1}.
\]
Combining this result with Lemma 2.2 we get 
\[
\frac{\alpha _l}{\alpha _{l-1}}=\frac{[n-l]_q}{[m-l]_q}\frac{q^{m-n}}{%
[m+n-2l+1]_q[m+n-2l+2]_q},l\geq 1.
\]
This immediately leads to the following expression of $\alpha _l$. 
\[
\alpha _l=q^{(m-n)l}\frac{[n]_q[n+1]_q[m+n]_q!}{[m]_q[m+1]_q}\prod_{i=1}^l%
\frac{[n-i]_q}{[m-i]_q}\prod_{i=1}^{2l}\frac 1{[m+n-2l+i]_q}.
\]
Comparing this expression with the previous one, after some simplification
we arrive at the following

{\bf Theorem 2.1.} Let q be an indeterminate and $m$ and $l$ be two positive
integers satisfying $m\geq l$. Then 
\[
\sum_{i=0}^{min\{l,m-l\}} (-1)^{i}q^{-i}\frac{[m-i]_{q}!}
{[i]_{q}![l-i]_{q}![m-l-i]_{q}}=q^{l(m-l)} 
\]
is a polynomial identity.

{\bf Corollary.} 
\[
\sum_{i=0}^{min\{l,m-l\}} (-1)^{i}C_{l}^{i}C_{m-i}^{l}=1.
\]

Proof. Rewrite Theorem 1 as 
\[
\sum_{i=0}^{min\{l,m-l\}}(-1)^{i}q^{-i} \left [%
\begin{array}{c}
l \\ 
i
\end{array}
\right ]_{q}\left [%
\begin{array}{c}
m-i \\ 
l
\end{array}
\right ]_{q}=q^{l(m-l)} 
\]
and take the classical limit.
\vspace{2mm}
{\bf 3. A New Basis of $V_{m}(x)\otimes V_{n}(y)$}
\vspace{1mm}

In this section we will establish a new basis of $V_{m}(x)\otimes V_{n}(y)$
under a certain condition. Without losing generality we assume $n\leq m$.

Let $j,l$ be two non-negative integers and $j\leq l\leq n$.We introduce the
notation 
\[
\phi _{l,j}=e_0^{l-j}f^j\Omega _0
\]
and define the sets 
\[
\Delta _l=\{\phi _{l,j}|j=0,1,\cdots ,l,l=0,1,\cdots ,n.\}
\]

Obviously, elements from different $\Delta _{l}$ are linearly independent as
they belong to different weight spaces. We will prove for each $l\in
\{0,1,\cdots,n\}$ $\Delta _{l}$ is a linearly independent set under some
condition. To this end, we will calculate explicitly the determinant of the
coefficient matrix of $\Delta _{l}$ with respect to the linearly independent
set $\{v_{i}\otimes w_{j}|i+j=l\}$.

Let 
\[
\phi _{l,j}=\sum_{i=0}^l\gamma _{l,j}^{i,l-i}v_i\otimes w_{l-i}.
\]
We denote by $(\Delta _l)$ the $l+1$ by $l+1$ coefficient matrix $(\gamma
_{l,j}^{i,l-i})$ whose row is marked by $j$ and column by $(i,l-i)$ and
denote by $|\Delta _l|$ the corresponding determinant. From the equations 
\[
\phi _{l+1,0}=e_0\phi _{l,0},\\ \phi _{l+1,j}=f\phi _{l,j-1},j=1,2,\cdots ,l,
\]
we get by direct calculation 
\begin{eqnarray*}
\gamma _{l+1,0}^{i,l+1-i} &=&yq^{-1}[l-i+1]_q\gamma
_{l,0}^{i,l-i}+xq^{-n+2l-2i+1}[i]_q\gamma _{l,0}^{i-1,l+1-i} \\
\gamma _{l+1,j}^{i,l-i+1} &=&q^{-m+2i}[l-i+1]_q\gamma
_{l,j-1}^{i,l-i}+[i]_q\gamma _{l,j-1}^{i-1,l-i+1},j=1,2,\cdots ,l+1.
\end{eqnarray*}

Multiply the column $(l+1,0)$ of $(\Delta _{l+1})$ by $(-q^{-m+2l})/[l+1]_q$
and add the result to the column $(l,0)$.Then multiply the column $(l,0)$ of
the resulted matrix by $(-q^{-m+2l-2})/[l]_q$ and add the result to the
column $(l-1,0)$.Continue this operation of canceling the "unwanted" terms
but keeping the value of the determinant untill we reach the first
column.This process leads to the following inductive formula: 
\[
\left| \Delta _{l+1}\right| =c[l+1]_q!\left| \Delta _l\right| ,
\]
where 
\[
c=(y-xq^{-m-n+2l})q^{-(l+1)}\sum_{i=0}^l(-1)^i\left[ 
\begin{array}{c}
l \\ 
i
\end{array}
\right] _qx^iy^{l-i}q^{i(l-1)}q^{-i(m+n)}.
\]
In deriving this formula,we have used the result: 
\[
\gamma _{l,0}^{l-i,i}=x^{l-i}y^i[l]_q!q^{-l}q^{(l-i)(i-n)},
\]
which can easily be obtained by straightforward calculation.

To get a neat expression for $c$ we need the following lemmas.

{\bf Lemma 3.1.} Let $l$ be a positive integer.For each $k\in \{l-1,l-3,\cdots,
l-1-2(l-1)\}$ 
\[
\sum_{i=0}^{l} (-1)^{i}\left [%
\begin{array}{c}
l \\ 
i
\end{array}
\right ]_{q} q^{ik}=0. 
\]

Proof.We  use induction method.When $l=1$ the formula is trivially
true.Using the formula 
\[
\left [%
\begin{array}{c}
l \\ 
i
\end{array}
\right ]_{q}=q^{-i}\left [ 
\begin{array}{c}
l-1 \\ 
i
\end{array}
\right ]_{q}+q^{l-i} \left [%
\begin{array}{c}
l-1 \\ 
i-1
\end{array}
\right ]_{q}, 
\]
we have 
\begin{eqnarray*}
&&\sum_{i=0}^{l} (-1)^{i}\left [%
\begin{array}{c}
l \\ 
i
\end{array}
\right ]_{q} q^{ik} \\
&&=\sum_{i=0}^{l-1} (-1)^{i}q^{(k-1)i}\left [%
\begin{array}{c}
l-1 \\ 
i
\end{array}
\right ]_{q}-q^{l+k-1}\sum_{j=0}^{l-1} (-1)^{j}q^{(k-1)j} \left [%
\begin{array}{c}
l-1 \\ 
j
\end{array}
\right ]_{q}.
\end{eqnarray*}
This implies the inductive step.Thus the lemma is proved.

{\bf Lemma 3.2.}Let $m,n$ and $l$ be integers,$l\geq 1$.Then for any complex
numbers $x,y$ 
\[
\sum_{i=0}^{l} (-1)^{i}\left [%
\begin{array}{c}
l \\ 
i
\end{array}
\right ]_{q} x^{i}y^{l-i}q^{i(l-1)}q^{-i(m+n)}\\=\prod_{j=0}^{l-1}
(y-xq^{-m-n+2j}). 
\]

Proof.Regard the left hand side as a polynomial of order $l$ in  indeterminate 
$y$. Then for each $j\in \{0,1,\cdots,l-1\}$ $xq^{-m-n+2j}$ is a root of
this polynomial thanks to Lemma 3.1. So the lemma follows.

With this lemma we are able to write down the determinant $|\Delta _{l+1}|$
in a neat form.

{\bf Proposition 3.1.}
\[
\left |\Delta _{l+1}\right |=\frac{[n]_{q}}{[m+1]_{q}}q^{-\frac{(l+1)(l+2)}
{2}}\prod_{j=1}^{l+1} [j]_{q}!\prod_{j=0}^{l} (y-xq^{-m-n+2j})^{l+1-j}.
\]

Proof.This is a direct consequence of the inductive formula and Lemma 3.2.

We are now finally prepared to prove the main result of this section.Let 
\[
\Delta =\bigcup_{l=0}^{n-1}\Delta _{l}
\bigcup_{l=0}^{n-1} f^{m+n-2l}\Delta _{l}\bigcup_{i=0}^{m-n} f^{i}
\Delta _{n}. 
\]
Here for a non-negative integer $i$ $f^{i}\Delta _{l}$ is defined to be the
set $\{f^{i}\phi _{l,j}|\phi _{l,j}\in \Delta _{l}\}$.

{\bf Theorem 3.1.}$\Delta$ is a basis of $V_{m}(x)\otimes V_{n}(y)$ if and only if
for \\each $j\in \{0,1,\cdots,n-1\}$, $y\neq xq^{-m-n+2j}$.

Proof.The "only if" part follows directly from Proposition 1. The same
proposition, together with the decomposition rule of $V_{m}\otimes V_{n}$
presented at the end of the first section, implies the "if" part.

Before concluding this paper we would like to state a dual form of Theorem
3.1. Let 
\begin{eqnarray*}
\varphi _{l,j}&=&f_{0}^{l-j}e^{j}\Phi _{0} \\
\Lambda _{l}&=&\{\varphi _ {l,j}|j=0,1,\cdots,l\},l=0,1,\cdots,n \\
\Lambda &=&\bigcup_{l=0}^{n-1}\Lambda _{l}
\bigcup_{l=0}^{n-1} e^{m+n-2l}\Lambda _{l}\bigcup_{i=0}^{m-n}
e^{i} \Lambda _{n}.
\end{eqnarray*}
We have the following

{\bf Theorem 3.2.}$\Lambda$ is a basis of $V_{m}(x)\otimes V_{n}(y)$ if and only
if for each $j\in \{0,1,\cdots,n-1\}$ $y\neq xq^{m+n-2j}$.

The proof of this theorem is parallel to that of the last theorem. We would
rather omit the details.
\newpage
Acknowledgement. The author thanks Dr. Omar Foda for the hospitality
extended to him during his stay in The University of Melbourne.\\

References\\

\begin{enumerate}
\item  V.Chari and A.N.Pressley,A Guide to Quantum Groups,Cambridge
University Press, Cambridge,1994.

\item  M.Jimbo,A q-analogue of U(gl(N+1)),Hecke algebra and the Yang-Baxter
equation.Lett.Math.Phys.11(1986)247.
\end{enumerate}

\end{document}